\documentclass[12pt]{amsart}

\usepackage{amsmath,amsthm,amsfonts,amssymb,eucal, hyperref}

\renewcommand {\a}{ \alpha }
\renewcommand{\b}{\beta}

\newcommand{\g}{\gamma}

\renewcommand{\d}{\delta}

\newcommand{\z}{\zeta}
\renewcommand{\t}{\theta}

\newcommand{\p}{\partial}

\newcommand{\Om}{\Omega}

\newcommand{\oq}{\ {\raise 7pt\hbox{${\scriptstyle\circ}$}}
\kern -7pt{
\hbox{$Q$}}}

\newcommand{\R}{ \mathbb R}

\newcommand {\bx}{\mathbf x}

\newcommand {\by}{\mathbf y}

\newcommand{\SPi}{{\sf{\Pi}}}
\newcommand{\SL}{{\sf{\Lambda}}}

\newcommand{\plainC}[1]{\textup{{\textsf{C}}}^{#1}}

\newcommand{\plainL}[1]{\textup{{\textsf{L}}}^{#1}}

\newcommand{\1}
{{\,\vrule depth3pt height9pt}{\vrule depth3pt height9pt} {\vrule
depth3pt height9pt}{\vrule depth3pt height9pt}\,}

 \DeclareMathOperator {\im }{{Im}}
\DeclareMathOperator {\re} {{Re}} \DeclareMathOperator {\dist}
{{dist}}

\newtheorem{thm}{Theorem}[section]
\newtheorem{cor}[thm]{Corollary}
\newtheorem{lem}[thm]{Lemma}
\newtheorem{prop}[thm]{Proposition}
\newtheorem{cond}[thm]{Condition}

\theoremstyle{definition}

\newtheorem{rem}[thm]{Remark}

\numberwithin{equation}{section}

\newcommand{\bee}{\begin{equation}}
\newcommand{\ene}{\end{equation}}
\newcommand{\bes}{\begin{split}}
\newcommand{\ens}{\end{split}}

\newcommand{\bet}{\begin{tm}}
\newcommand{\ent}{\end{tm}}
\newcommand{\bel}{\begin{lm}}
\newcommand{\enl}{\end{lm}}
\newcommand{\bec}{\begin{cor}}
\newcommand{\enc}{\end{cor}}
\newcommand{\bep}{\begin{pr}}
\newcommand{\enp}{\end{pr}}
\newcommand{\ber}{\begin{rem}}
\newcommand{\enr}{\end{rem}}

\makeatletter
\def\square{\RIfM@\bgroup\else$\bgroup\aftergroup$\fi
  \vcenter{\hrule\hbox{\vrule\@height.6em\kern.6em\vrule}\hrule}\egroup}
\makeatother

\begin{document}

\hoffset -4pc

\title[Hardy inequalities]
{Hardy inequalities for simply connected planar domains}
\author[A. Laptev]{Ari Laptev}
\author[A. V. Sobolev]{Alexander V. Sobolev}

\address{Department of Mathematics\\
Royal Institute of Technology\\
SE-100 44 Stockholm\\
Sweden} \email{laptev@math.kth.se}
\address{School of Mathematics\\
University of Birmingham\\
Edgbaston, Birmingham\\ B15 2TT, UK}
\email{A.Sobolev@bham.ac.uk}

\date{\today}

\begin{abstract}
In 1986 A. Ancona showed, using the Koebe one-quarter Theorem, that
for a simply-connected planar domain the constant in the Hardy
inequality with the distance to the boundary  is greater than or equal
to $1/16$. In this paper we consider classes of domains for which
there is a stronger version of the Koebe Theorem. This
implies better estimates for the constant appearing in the Hardy
inequality.
\end{abstract}

\maketitle
\vskip 0.5cm

\section{Main result and discussion}

Let $\Om$ be a domain in $\mathbb R^2$ and let $\Om^c =\mathbb R^2\setminus \Om$ be its complement.
 For any function $u\in
\plainC{1}_0(\Om)$ we have:
\begin{equation}\label{hardy:eq}
\int_{\Om} |\nabla u|^2 d\bx \ge r^2 \int_{\Om} \frac{|u|^2}{\d(\bx)^2}
d\bx,\quad \d(\bx) = \inf_{\by\in\Om^c}|\by-\bx|,
\end{equation}
e.g. E.B. Davies  \cite{D1}, \cite{D2}, \cite{D3} and V.G. Maz'ya \cite{M}.
It is well known that for convex domains  $r = 1/2$ and it is sharp, see e.g.
\cite{D1}.
However, the sharp constant for non-convex domains
is unknown,  although for arbitrary planar simply-connected domains
A. Ancona \cite{A} proved \eqref{hardy:eq} with $r=1/4$.
Some specific examples of non-convex domains were considered
in \cite{D3}  (see also J. Tidblom \cite{T2}).
For example, it was found that if $\Om={\mathbb R}^2\setminus {\mathbb R}_+$,
${\mathbb R}_+ = [0,\infty)$, then $r^2= 0.20538...$.

Our objective is to obtain the Hardy inequality for
simply-connected non-convex domains $\Om\subset\R^2$, whose degree of non-convexity
can be "quantified". We introduce two possible "measures" of
non-convexity.

Let $\SL\subset\mathbb C$ be a simply-connected domain such that $0\subset\p\SL$.
Denote by
$\SL(w, \phi) = e^{i\phi}\SL+w$
the transformation of $\SL$ by rotation by angle $\phi\in (-\pi, \pi]$
in the positive direction and translation by $w\in\mathbb C$:
\begin{equation}\label{trans:eq}
\SL(w, \phi) = \{ z\in\mathbb C: e^{-i\phi}(z-w)\in \SL\}.
\end{equation}
Denote by $K_{\t}\subset \mathbb C$, \ $\t\in [0, \pi]$ the sector
\begin{equation}\label{sector:eq}
K_{\t} = \{z\in \mathbb C:  |\arg z|< \t\}.
\end{equation}
In words, this is an open sector symmetric with respect to the
real axis, with the angle $2\t$ at the vertex.
Here and below we always assume that $\textup{arg}\
\z\in (-\pi, \pi]$ for all $\z\in\mathbb C$. Our first assumption on
the domain $\Om$ is the following

\begin{cond}\label{cone:cond} There exists a number $\t\in [0, \pi]$
such that for each $w\in \Om^c$ one can find a $\phi = \phi_w\in
(-\pi, \pi]$ such that
\begin{equation*}
\Om\subset K_{\t}(w, \phi_w).
\end{equation*}
\end{cond}

Very loosely speaking, this means that the domain $\Om$ satisfies the
exterior cone condition. The difference is of course that the cone
is now supposed to be infinite. Because of this, Condition \ref{cone:cond} is equivalent
to itself if stated
for the boundary points $w\in\p\Om$ only.

Note also that if Condition
\ref{cone:cond} is satisfied for some $\t$, then automatically
$\t\ge\pi/2$, and the equality $\t = \pi/2$ holds for convex domains.

\begin{thm}\label{main1:thm}
Suppose that the domain $\Om\subset\R^2, \Om\not = \R^2$
satisfies Condition \ref{cone:cond} with some $\t\in [\pi/2, \pi]$.
Then for any $u\in \plainC{1}_0(\Om)$ the Hardy inequality
\eqref{hardy:eq} holds with
\begin{equation}\label{r1:eq}
r = \frac{\pi}{4\t}.
\end{equation}
\end{thm}

It is clear that the constant $r$
runs from $1/4$ to $1/2$ when $\t$ varies from $\pi$ to $\pi/2$. For the
domain $\Om = K_{\t}$ Theorem \ref{main1:thm} does not give the best known result,
found in \cite{D3}, saying that the value of $r$ remains equal to $1/2$ for the range
$\t \in [0, \t_0]$ where $\t_0 \approx 2.428$, which is considerably greater than $\pi/2$.

To describe another way to characterize the non-convexity, for
$a >0$ and $\t\in [0, \pi)$, introduce the domains
\begin{equation}\label{pan:eq}
\tilde D_{a} = \{z\in\mathbb C:  |z|> a\  \& \ \
|\arg z|\not = \pi\},\ \   D_{a, \t} = \tilde D_{a}(-a e^{i\t}, 0).
\end{equation}
The domain $\tilde D_{a}$ is the exterior of the disk of radius $a$
centered at the origin with an infinite cut along the negative real semi-axis.

\begin{cond} \label{outside:cond}
There exist numbers $a >0$ and $\t_0\in [0, \pi)$
such that for any $w\in \p\Om$ one can
find a $\phi = \phi_w\in (-\pi, \pi]$ and $\t\in [0, \t_0]$ such that
\begin{equation*}
\Om\subset  D_{a, \t}(w, \phi_w).
\end{equation*}
\end{cond}

Note that any domain satisfying Conditions \ref{cone:cond} or
\ref{outside:cond}, is automatically simply-connected.

The following Theorem applies to the domains with a finite in-radius
\begin{equation*}
\d_{\textup{\tiny in}} = \sup_{z\in\Om} \d(z).
\end{equation*}

\begin{thm}\label{main2:thm}
Let $\Om\subset\R^2$, $\Om\not=\R^2$ be a domain such that
$\d_{\textup{ in}} <\infty$. Suppose that $\Om$ satisfies Condition
\ref{outside:cond} with some $\t_0\in [0, \pi)$ and that
$$
2\d_{\textup{\tiny in}}\le  R_0(a),\ \ R_0(a) = \frac{a}{2(2^{1/2} |\tan(\t_0/2)| + 1)}.
$$
Then the Hardy
inequality \eqref{hardy:eq} holds with
\begin{equation}\label{r2:eq}
r = \frac{1}{2}
\biggl[
1 - 4 \bigl(2^{1/2}|\tan(\t_0/2)| + 1\bigr)\frac{\d_{\textup{in}}}{a}
\biggr].
%
%
\end{equation}
\end{thm}

A natural example of a domain to apply the above
theorem, is the following horseshoe-shaped domain
$$
\SL = \{z\in\mathbb C: \rho<|z|< \rho + \d,\ |\arg z| < \psi,\
\re z > \rho\cos\psi\},\ \ \psi\in (0, \pi),
$$
with $\rho, \d >0$.
Simple geometric considerations show that
this domain satisfies Condition
\ref{outside:cond} with $a = \rho$ and
\begin{equation*}
\t_0 = \begin{cases}
0,\ \ \psi\le \pi/2,\\
2\psi -\pi,\ \ \psi > \pi/2.
\end{cases}
\end{equation*}
Assuming that $\d \rho^{-1}$ is small, so that $\d_{\textup{in}} = \d$,
we deduce from Theorem \ref{main2:thm} that the Hardy inequality holds
with a constant $r$, which gets close to $1/2$ as $\d \rho^{-1}\to 0$.
On the other hand, if $\d_{\textup{in}}\rho^{-1}$ is large, one could apply Theorem
\ref{main1:thm},
noticing that $\SL$ satisfies Condition \ref{cone:cond}
with $\t = (\pi+\psi)/2$, which gives the
Hardy inequality with constant
$$
r = \frac{1}{2\bigl(1+\psi\pi^{-1}\bigr)}.
$$
which is obviously independent of $\d_{\textup{in}}$ or $\rho$.

Let us mention briefly some other recent results for convex domains,
concerning the Hardy inequality
with a remainder term. In the paper \cite{BM} H.Brezis and M. Marcus
showed that if $\Om\in{\mathbb R}^d$, $d\ge 2$,
then the inequality  could be improved to include the $\plainL2$-norm:
\begin{equation}\label{hardy:eq:corr}
\int_{\Om} |\nabla u|^2 d\bx \ge \frac{1}{4}\int_{\Om} \frac{|u|^2}{\d(\bx)^2}
+ C(\Om) \,\int_{\Om} |u|^2 d\bx,
\end{equation}
where the constant $C(\Om)>0$ depends on the diameter of $\Om$.
They also conjectured that $C(\Om)$ should depend on
the Lebesgue measure of $\Om$. This conjecture was justified in \cite{HOHOL}
and later generalised to $\plainL{p}$-type inequalities  in \cite{T1}.
Later S. Filippas,  V.G. Maz'ya and A. Tertikas \cite{FMT} (see also
F.G. Avkhadiev \cite{A}) obtained
for $C(\Om)$ an estimate in terms of the in-radius $\d_{\textup{in}}$.

\section{A version of the Koebe Theorem}

A. Ancona has pointed out in \cite{A} (page 278)
that the Hardy inequality for simply-connected planar domains can be obtained from the famous Koebe
one-quarter Theorem. Let $f$ be a
conformal mapping (i.e. analytic univalent) defined on the unit disk
$\mathbb D = \{z\in\mathbb C: |z| < 1\}$, normalized by the
condition $f(0) = 0$, $f'(0) = 1$.
 Denote by  $\Om$ the image of the disk under the function $f$,
i.e. $\Om = f(\mathbb D)$, and set
$$
\d(\z) = \dist\{\z, \p\Om\} = \inf_{w\notin\Om}|w-\z|
$$
to be the distance from the point $\z\in\Om$ to the boundary
$\p\Om$. The classical Koebe  one-quarter Theorem tells us that
\begin{equation*}
\d(0) \ge r,
\end{equation*}
with $r = 1/4$. On the other hand,  if the domain $\Om$ is convex,
then it is known that $r = 1/2$, see e.g. P.L.Duren \cite{Duren}, Theorem 2.15.
Without the  normalization conditions $f(0) = 0$, $f'(0) = 1$ the
above estimate can be rewritten as follows:
\begin{equation}\label{koebe:eq}
\d(f(0))\ge r |f'(0)|.
\end{equation}
For any simply-connected domain $\Om\subset \mathbb C,\ \Om\not = \mathbb C$ we denote
by $\mathbb A(\Om)$ the class of all conformal maps such that
$f(\mathbb D) = \Om$.

Our proof of the main Theorems \ref{main1:thm} and \ref{main2:thm} relies on a
version of the Koebe theorem, in which the constant $r$ assumes values in the interval
$[1/4, 1/2]$. We begin with a general statement which deduces the required Koebe-type result
by comparing the domain $\Om$ with some suitable "reference" domain.
Let $\SL\subset \mathbb C, \SL \not = \mathbb C$
be a simply-connected domain such that $0\subset\p\SL$, and let
$g$ be a conformal function which maps $\SL$ onto the complex
plane with a cut along the negative semi-axis, i.e. onto
\begin{equation*}
\SPi = \mathbb C\setminus\{z\in\mathbb C: \im z = 0, \re z \le 0\},
\end{equation*}
such that $g(0) = 0$. We call $\SL$  a \textsl{standard domain} and $g$ -
a \textsl{comformal map associated} with the standard domain $\SL$.

\begin{lem}\label{gprime:lem}
Let $w\in\p\Om$ and suppose that for some standard
domain $\SL$ the inclusion
\begin{equation}\label{standard:eq}
\Om\subset\SL(w, \phi)
\end{equation}
holds with some $\phi\in (- \pi, \pi]$.
Let $g$ be a conformal map associated with $\SL$, and suppose that there are numbers
$M\in (0, \infty)$
and $R_0\in (0, \infty]$ such that
for all $R\in (0, R_0)$
\begin{equation}\label{gprime:eq}
\left|\frac{g'(z)}{g(z)}\right|\ge \frac{\b}{|z|}
\end{equation}
for $z: 0< |z| \le R, z\in\SL$ with some constant $\b = \b(R) \in (0, M]$.
Then for any $f\in \mathbb A(\Om)$ satisfying the
condition $M|f'(0)| < 4 R_0$, the inequality
\begin{equation}\label{koebe2:eq}
|f(0) - w|\ge \frac{\b(R_1)}{4}|f'(0)|,\ \  R_1 = \frac{M|f'(0)|}{4},
\end{equation}
holds.
\end{lem}

\begin{proof}
Since $\Om\subset\SL(w, \phi)$, the function
\begin{equation*}
h(z) = g\bigl(e^{-i\phi}(f(z)-w)\bigr)
\end{equation*}
is conformal on $\mathbb D$.
Since $0\notin \SPi$, by the classical Koebe Theorem,
\begin{equation*}
|h(0)|\ge \frac{1}{4} |h'(0)|,
\end{equation*}
so that
\begin{equation*}
|g\bigl(e^{-i\phi}(f(0)-w)\bigr)|
\ge \frac{1}{4}|g'\bigl(e^{-i\phi}(f(0)-w)\bigr)| |f'(0)|.
\end{equation*}
If $|f(0)-w|\ge M |f'(0)|/4$,
then there is nothing to prove, so we assume that
$|f(0) - w|\le M|f'(0)|/4$.
Then by the assumption \eqref{gprime:eq} we get
\begin{equation*}
\frac{1}{\b\bigl(M|f'(0)|/4\bigr)}|f(0) - w|\ge \frac{1}{4}|f'(0)|,
\end{equation*}
which leads to the required estimate.
\end{proof}

\begin{cor}\label{standard:cor}
Let $\Om\subset\mathbb C$ be a domain.
Suppose that for any $w\in \p\Om$ there is a standard domain
$\SL = \SL_w$, such that the inclusion \eqref{standard:eq} holds with some $\phi = \phi_w$,
and that the associated conformal maps $g = g_w$ satisfy \eqref{gprime:eq}
for all $0 <|z|\le R, z\in\SL_w$, for all $R < R_0$
with some $\b(R) = \b_w(R)\in (0, M]$ where $M\in (0, \infty)$ and $R_0\in (0, \infty]$
are independent of $w$.

If $R_0 <\infty$, then under the condition $M\d_{\textup{in}} < R_0$ the
estimate \eqref{koebe:eq} holds for all $f\in \mathbb A(\Om)$ with
$$
r = \frac{1}{4}\ \inf_{w\in\p\Om}\b_w(R'),\ R' = M \d_{\textup{in}}.
$$

If $R_0 = \infty$, then the estimate \eqref{koebe:eq} holds for all $f\in\mathbb A(\Om)$
with
$$
r = \frac{1}{4}\ \inf_{w\in\p\Om}\inf_{R>0}\b_w(R).
$$
\end{cor}

Observe that under the conditions of this corollary, the domain $\Om$
is automatically simply-connected and $\Om\not = \mathbb C$.

\begin{proof}
In the case $R_0 = \infty$ the result immediately follows
from Lemma \ref{gprime:lem}.

Assume that $R_0 <\infty$.
By the classical Koebe Theorem $|f'(0)|\le 4\d(f(0))\le 4\d_{\textup{in}}$, so that
by Lemma \ref{gprime:lem}, for each $w\in\p\Om$ we have the estimate
\eqref{koebe2:eq}. Since $R_1 \le R'$ and $\b_w(\ \cdot\ )$ is a decreasing function,
the required inequality \eqref{koebe:eq} follows.
\end{proof}

Now we apply the above results in the cases of standard domains $K_{\t}$ and
$D_{a, \t}$, see \eqref{sector:eq} and \eqref{pan:eq} for definitions.


\begin{thm}\label{koebe:thm}
Suppose that $\Om$ satisfies Condition \ref{cone:cond}
with some $\t\in [\pi/2, \pi]$.
Then for any $f\in\mathbb A(\Om)$ the inequality
\eqref{koebe:eq} holds with $r$ given by \eqref{r1:eq}.
\end{thm}

\begin{proof}
Due to Condition \ref{outside:cond}, for each $w\in\p\Om$ we have
$\Om\subset K_{\t}(w, \phi)$ with some $\phi\in(-\pi, \pi]$.
Clearly, the domain $K_{\t}$ is standard and the function
$$
g(z) = z^{\a},\ \a = \frac{\pi}{\t}
$$
is a conformal map associated with $K_{\t}$.
One immediately obtains:
$$
\frac{g'(z)}{g(z)} = \a \frac{z^{\a-1}}{z^\a} = \frac{\a}{z}, z\in \SL,
$$
so that the conditions of Corollary
\ref{standard:cor} hold with the constant $\b=\a$
and $R_0 = \infty$.
Now Corollary \ref{standard:cor} leads to the proclaimed result.
\end{proof}

Note that for convex domains the angle $\t$ is $\pi/2$, and hence we
recover the known result $r = 1/2$. Actually, the proof of Lemma \ref{gprime:lem}
is modelled on that for convex domains, which is  featured in
\cite{Duren}, Theorem 2.15.

Let us prove a similar result for Condition \ref{outside:cond}:

\begin{thm}\label{pan:thm}
Suppose that $\Om$ satisfies Condition
\ref{outside:cond} with some $a >0$,\ $\t_0\in [0, \pi)$,
and that $\d_{\textup{in}} < \infty$,
\begin{equation}\label{din:eq}
2\d_{\textup{in}}\le R_0(a), \ \ R_0(a) = \frac{a}{2\bigl(
2^{1/2}|\tan(\t_0/2)| + 1
\bigr)}.
\end{equation}
Then for any $f\in\mathbb A(\Om)$ the inequality
\eqref{koebe:eq} holds with $r$ given by \eqref{r2:eq}.
\end{thm}

\begin{proof}
Due to Condition \ref{outside:cond}, for each $w\in\p\Om$ we have
$\Om\subset D_{a, \t}(w, \phi)$ with some $\t\in [0, \t_0]$ and $\phi\in(-\pi, \pi]$.
Clearly, the domain $D_{a, \t}$ is standard and the function
$$
g(z) = h^2(za^{-1}),\ \
h(\z) =  \sqrt{\z+ e^{i\t}} - \frac{1}{\sqrt{\z+ e^{i\t}}}
- 2i b,\ \ b = \sin \frac{\t}{2},
$$
is a conformal map associated with $D_{a, \t}$.
Write:
\begin{align*}
h(\z) = &\ \frac{\z+e^{i\t} - 1
- 2ib \sqrt{ \z+  e^{i\t} }}{\sqrt{ \z+ e^{i\t} }}\\[0.3cm]
= &\ \frac{\z+  e^{i\t} - 1
- 2ib \sqrt{ \z+  e^{i\t} }}{\sqrt{ \z+  e^{i\t} }}\\[0.3cm]
= &\ \frac{\z+ 2i  b e^{i\t/2}
- 2ib \sqrt{ \z+  e^{i\t} }}{\sqrt{ \z+  e^{i\t}}}\\[0.3cm]
= &\ \frac{\psi(\z)}{\sqrt{\z+e^{i\t}}}.
\end{align*}
A direct calculation shows that
\begin{equation*}
\frac{g'(z)}{g(z)} = \frac{2\psi'(z/a)}{a\psi(z/a)} - \frac{1}{z+a e^{i\t}}.
\end{equation*}
Let us investigate the function $\psi$ is more detail. Assume that $|\z|\le 1/2$.
 Rewrite:
\begin{align*}
\psi(\z) = &\ \z + 2ib e^{i\t/2} - 2ib e^{i\t/2}\sqrt{1 + \z e^{-i\t/2}}\\[0.3cm]
= &\ \z + 2ib e^{i\t/2} - 2ibe^{i\t/2}\biggl(1+ \frac{1}{2}\z e^{-i\t}+ \z\g_1\biggr)\\[0.3cm]
= &\ (1-ibe^{-i\t/2}) \z - 2 ib e^{i\t/2} \z\g_1
= e^{-i\t/2} \cos \frac{\t}{2}\  \z - 2 ib e^{i\t/2} \z\g_1,
\end{align*}
where
$$
|\g_1(\z)|\le 2^{-3/2} |\z|,\ |\z|\le 1/2.
$$
Let's look at the derivative:
\begin{align*}
\psi'(\z) = &\ 1 - \frac{ib}{\sqrt{\z+e^{i\t}}}\\
= &\ 1 - e^{-i\t/2}\frac{ib}{\sqrt{1+\z e^{-i\t}}}
= 1 - ib e^{-i\t/2}(1+\g_2)\\
= &\  e^{-i\t/2} \cos \frac{\t}{2} - ib e^{-i\t/2}\g_2,
\end{align*}
where
$$
|\g_2(\z)|\le 2^{1/2}|\z|, \ |\z|\le 1/2.
$$
Therefore
\begin{align*}
\frac{g'(z)}{g(z)} = &\ \frac{2e^{-i\t/2} \cos \frac{\t}{2} - 2ib e^{-i\t/2}\g_2(z/a)}
{e^{-i\t/2} \cos \frac{\t}{2}\  z - 2 i  b e^{i\t/2} z\g_1(z/a)} - \frac{1}{z+a e^{i\t}}\\[0.3cm]
= &\ \frac{1}{z}\biggl[
\frac{2 - 2i\tan(\t/2) \g_2(z/a)}{1 - 2i \tan (\t/2) e^{i\t}\g_1(z/a)}
- \frac{z}{a} \frac{1}{z/a + e^{i\t}}
\biggr]
\end{align*}
A simple calculation shows that
\begin{align*}
\left|\frac{g'(z)}{g(z)}\right| \ge
 &\ \frac{2}{|z|}\biggl[
\frac{1 - |\tan(\t/2)| |\g_2(z/a)|}{1 + 2 |\tan (\t/2)| |\g_1(z/a)|}
- 2\frac{|z|}{a}
\biggr]\\[0.3cm]
\ge &\ \frac{2}{|z|}\biggl[
1 - 2 \bigl(2^{1/2}|\tan(\t/2)| + 1\bigr)\frac{|z|}{a}
\biggr], \ |z|/a\le 1/2.
\end{align*}
Therefore the condition \eqref{gprime:eq} is satisfied for all $0<|z|\le R$,
$$
R < R_0 = \frac{a}{2\bigl(
2^{1/2}|\tan(\t/2)| + 1
\bigr)},
$$
with
\begin{equation}\label{beta:eq}
\b(R) = \b_{\t}(R) = 2\biggl[
1 - 2 \bigl(2^{1/2}|\tan(\t/2)| + 1\bigr)\frac{R}{a}
\biggr].
\end{equation}
Note that $\b\le M$ with $M=2$ and $\b_{\t}\ge \b_{\t_0}$, so that by Corollary \ref{standard:cor}
the estimate \eqref{koebe:eq} holds with the  $r$ given by \eqref{r2:eq}.
\end{proof}

\section{Proof of Theorems \ref{main1:thm}, \ref{main2:thm}}

As soon as the Koebe Theorem \eqref{koebe:eq} is established, our
proof of the Hardy inequality follows that by A.Ancona \cite{A}.
Namely, our starting point is the inequality for the
half-plane, which is an immediate consequence of
the classical Hardy inequality in one dimension.
Below we use the usual notation $z = x+ iy$, $x, y\in \R$.

\begin{prop}\label{hplane:prop} Let $\mathbb C_+ = \{z\in\mathbb C: \re z >0\}$.
For any $u\in \plainC{1}_0(\mathbb C_+)$ one has
\begin{equation*}
\int_{\mathbb C_+} |\nabla u|^2 dx dy \ge \frac{1}{4} \int_{\mathbb
C_+} \frac{|u|^2}{x^2} dx dy.
\end{equation*}
\end{prop}

Theorems \ref{main1:thm} and \ref{main2:thm} immediately follow
from Proposition \ref{hplane:prop} with the help of the following
conditional result:

\begin{thm}\label{hkoebe:thm}
Let $\Om\subset\mathbb C$, $\Om\not=\mathbb C$ be a simply connected domain.
Suppose that for all $g\in \mathbb A(\Om)$ the inequality
\eqref{koebe:eq} holds with some $r\in [1/4, 1/2]$. Then for any
conformal mapping $f:\mathbb C_+\to \Om$ the following version of
the Koebe Theorem holds:
\begin{equation}\label{dist:eq}
\d(f(z))\ge 2rx|f'(z)|.
\end{equation}
\end{thm}

For $r = 1/4$ the estimate \eqref{dist:eq} can be found in \cite{A}.
For the reader's convenience we provide a proof of \eqref{dist:eq}.

\begin{proof}
For a conformal mapping  $f: \mathbb C_+\to \Om$ and arbitrary $z\in
\mathbb C_+$ we define
\begin{equation*}
g(w) = g_z(w) = f\bigl(h(w)\bigr),\ h(w) = \frac{\overline z w +
z}{1-w} ,
\end{equation*}
where $w\in \mathbb D$. It is clear that for each fixed $z\in\mathbb
C_+$ the function $h$ maps $\mathbb D$ onto $\mathbb C_+$ and $h(0)
= z, \ g(0) = f(z)$. The derivative is
\begin{equation*}
g'(w) = \frac{z+\overline z}{(1-w)^2} f'\bigl(h(w)\bigr),
\end{equation*}
so that
\begin{equation*}
g'(0) = 2 x f'(z).
\end{equation*}
Therefore, the Koebe theorem \eqref{koebe:eq} implies that
\begin{equation*}
\d\bigl(f(z)\bigr) = \d\bigl(g(0)\bigr)\ge r |g'(0)| = 2r x |f'(z)|,
\end{equation*}
as required.
\end{proof}

\begin{proof}[Proof of Theorems \ref{main1:thm}, \ref{main2:thm}]
According to Theorems \ref{koebe:thm} or \ref{pan:thm} the Koebe
Theorem for functions $f\in \mathbb A(\Om)$ holds with the values of
$r$ given by \eqref{r1:eq} or \eqref{r2:eq} respectively.

Let $f:\mathbb C_+\to \Om$ be a conformal map. Remembering that
conformal maps preserve the Dirichlet integral, from Proposition
\ref{hplane:prop} we get for any $u\in \plainC{1}_0(\Om)$
\begin{align*}
\int_{\Om}|\nabla u|^2 d\bx = &\ \int_{\mathbb C_+} |\nabla(u\circ
f)|^2 d x dy
\ge \frac{1}{4} \int_{\mathbb C_+}  \frac{|(u\circ f)|^2}{x^2} d\bx \\[0.3cm]
=  &\ r^2 \int_{\mathbb C_+}  \frac{|(u\circ f)|^2}{ (2rx)^2
|f'(z)|^2}  |f'(z)|^2  d\bx \ge  r^2\int_{\Om}
\frac{|u|^2}{\d(\bx)^2} d\bx.
\end{align*}
At the last step we have used \eqref{dist:eq}.

Now Theorems \ref{main1:thm}, \ref{main2:thm} follow.
\end{proof}

\medskip
\noindent
{\it Acknowledgements.}
The authors are grateful to partial support by the SPECT ESF European programme.
A. Laptev would like to thank the University of Birmingham for its hospitality.
A. Sobolev was partially supported by G\"oran Gustafsson Foundation.

\bibliographystyle{amsplain}

\providecommand{\bysame} {\leavevmode\hbox
to3em{\hrulefill}\thinspace}

\end{document}